\documentclass[11pt,reqno]{amsart}
\usepackage{latexsym}
\newtheorem{lem}{Lemma}
\newtheorem{thm}{Theorem}
\newtheorem{prop}[thm]{Proposition}
\newtheorem{cor}[thm]{Corollary}
\theoremstyle{remark}
\newtheorem{rem}{Remark}
\theoremstyle{definition}
\newtheorem{defi}{Definition}
\newcommand{\R}{\mathbb{R}}
\newcommand{\C}{\mathbb{C}}
\newcommand{\K}{\mathbb{K}}
\newcommand{\N}{\mathbb{N}}

\newcommand{\vpol}[1][*]{\mathit{Vect}_{{#1}}} 
\newcommand{\euler}{\mathcal{E}}
\newcommand{\ad}{\mathit{ad}}
\newcommand{\del}{\partial}

\newcommand{\alg}[2]{\mathit{{#1}}({#2})}

\newcommand{\descc}[1][]{\mathit{c}^{#1}}
\newcommand{\n}[1][]{\mathcal{N}^{#1}}
\newcommand{\T}[1][*]{\mathcal{T}_{{#1}}}

\newcommand{\id}{\mathrm{id}}

\makeatletter
\def\@cite#1#2{{%
  \m@th\upshape\mdseries Ref.~{#1\if@tempswa, #2\fi}}}
\@ifundefined{cite }{%
  \expandafter\let\csname cite \endcsname\cite
  \edef\cite{\@nx\protect\@xp\@nx\csname cite \endcsname}%
}{}
\makeatother

\title[Maximal subalgebras]{Maximal subalgebras of  vector 
fields for equivariant quantizations}
\author{F.  Boniver and P.  Mathonet}

\begin{document}
\begin{abstract}
    The elaboration of new quantization methods has recently 
    developed the interest in the study of subalgebras of the Lie algebra of polynomial 
    vector fields over a Euclidean space.  In this framework, these subalgebras define 
    maximal equivariance conditions that one can impose on a linear 
    bijection between observables that are polynomial in the momenta 
    and differential operators.
    Here, we determine which finite dimensional graded Lie subalgebras 
    are maximal.
    In order to characterize these, we make use of results of Guillemin, Singer and 
    Sternberg  and Kobayashi and Nagano. 
\end{abstract}
    \maketitle
PACS classification numbers: 03.65.F, 02.10.Vr, 02.10.Sp
    \newpage
\section{Introduction}
Our interest in the present study comes from recent works about 
new equivariant quantizations~(\cite{dlo,lo}).

One can define \emph{quantization maps}
as linear bijections~$\mathcal{Q}$ 
from the space $\mathit{Pol}(T^* M)$ of functions 
on the cotangent bundle of a smooth 
manifold $M$, that are polynomial on the fibre, to a space 
$D_{\lambda}(M)$ of 
differential operators acting on tensor densities of weight $\lambda$ 
over $M$.

It is known that a quantization map $\mathcal{Q}$ cannot be equivariant 
  with respect to all diffeomorphisms of $M$.
From the infinitesimal point of view, this means that such a map does 
not commute with the action on these spaces of
the Lie algebra $\vpol[](M)$
 of vector fields over $M$.
In other words, differential operators and polynomials are 
inequivalent  modules of  $\vpol[](M)$.

However, when $M$ is endowed with an additional structure, some
 particular subalgebras of $\vpol[](M)$ naturally deserve 
 consideration, because they are made up of 
infinitesimal transformations preserving the structure.

The authors of~\cite{dlo,lo} considered the case of infinitesimal 
projective or conformal transformations of $M$.
In suitable charts, these can be realized in 
polynomial vector fields over a Euclidean space.  
For instance, if $M$ is endowed with a projective 
structure (i.e. $M$ is locally identified with a real projective 
space, say of dimension $n$) then in appropriate charts, the Lie algebra 
of infinitesimal projective transformations -- isomorphic to 
$\alg{sl}{n+1,\R}$ -- is generated by the vector fields
\begin{equation}\label{sln}
\frac{d}{d x^j}, \; x^j\frac{d}{dx^k},\; x^j\sum_{l=1}^n x^l 
\frac{d}{dx^l},\quad \forall j,k\leq n.
\end{equation}
 
 In this setting, 
those conformal and projective subalgebras share the property of being maximal 
in the algebra of polynomial vector fields: they are not contained in 
any larger proper subalgebra.  
The reader may refer to~\cite{lo, bolec} for proofs.

Now, it was proved in~\cite{dlo,lo} that one could
 construct a  quantization map  equivariant with 
respect to those subalgebras.  This 
quantization is unique up to normalization.

In this framework, our concern in the present paper is to determine all finite 
dimensional graded 
subalgebras of polynomial vector fields over a given Euclidean space 
that are maximal.

Independently of 
quantization purposes, other maximality conditions have also been studied. 

In~\cite{kan}, Kantor classified \emph{irreducible transitively 
differential groups}. This notion
gives rise, from the Lie algebraic point of view, to the class of finite 
dimensional graded Lie subalgebras of polynomial vector fields 
containing all constant vector fields.
The author then seeks for irreducible~(see \cite[p.~1405]{kan} or 
below) subalgebras being maximal in 
this class.

Another more recent study is that of Post~(\cite{post}).  In this 
paper, a stronger grading
requirement is imposed  in order to define a class 
of finite dimensional Lie algebras containing all constant vector 
fields.
All maximal subalgebras of this class are then identified.

We point out two differences between the maximality conditions examined 
here and in these studies. 

On the one hand, we impose fewer conditions on the subalgebras we 
consider, keeping only the requirements for a subalgebra to be graded 
and finite dimensional.
On the other hand, the maximality property is not investigated inside 
a particular class of subalgebras, but in the general class of all 
subalgebras of polynomial vector fields.

Before giving  our main result and a brief description of the tools we shall use,  
let us fix some notations.

 Throughout this note, we assume that $E$ is an $n$-dimensional 
vector space over $\K$, which is taken to be $\R$ or $\C$.  We shall 
deal with polynomial vector fields over $E$.

We denote by $\vpol (E)$ the space of these vector fields, i.e. the 
space of polynomial maps from $E$ to $E$.  It is worth noticing that 
the vector fields considered when $E$ is complex are thus holomorphic.
Let $\{e_{j},j=1,\ldots,n\}$ be a basis of $E$.
Assume that $X,Y\in \vpol(E)$ are written
$X =\sum_{j=1}^n X^j e_{j}$ and $Y= \sum_{j=1}^n Y^j e_{j}$.
We denote as usual by $[X,Y]$ the Lie bracket
\[
\sum_{j,k} X^j\del_{j} Y^k e_{k}- Y^{j}\del_{j} X^k e_{k},
\]
where $\del_{j}$ represents the derivation $\frac{d}{dx^j}$ along 
the $j$-th axis.  For the sake of convenience, we shall also use this 
notation to designate the $j$th vector of a basis of $E$.
We denote by $\ad(X)$ the map $Y \mapsto [X, Y]$.

We name \emph{Euler vector field} the
identity transformation of $E$.  In a basis $\{\del_{j}\}$, it reads
\[
\euler(x)=\sum x^j\del_{j}.
\]
It defines a natural grading on $\vpol(E)$~:
\[
\vpol(E) = \bigoplus_{p\geq -1} \vpol[p](E)
\]
where $\vpol[p](E)$ denotes the space of eigenvectors of
$\ad(\euler)$ associated with the eigenvalue $p$, i. e. vector fields
with homogeneous coefficients of degree $p+1$.
 
We are interested in these \emph{graded} subalgebras $L$ :
\[
L=\bigoplus_{-1\leq p\leq r}L_{p} \quad\text{ with }
L_{p}=\vpol[p](E)\cap L
\]
which are maximal in $\vpol (E)$.  As mentioned above, the 
notion of maximality has been used in various senses.  Therefore, it is  
worth emphasizing the following definition.
\begin{defi}
    A subalgebra $L$ of $\vpol(E)$ is \emph{maximal} if 
    \begin{eqnarray*}
    L\subset L'&\Rightarrow& L'=L \\
    &&\text{ or }L'=\vpol(E)
    \end{eqnarray*}
    whenever $L'$ is a subalgebra of $\vpol(E)$.
\end{defi}

\section{Main result}
\begin{defi}[\protect{see for instance \cite[p.~682]{koba3}}]
    A graded subalgebra $L=L_{-1}\oplus\cdots\oplus L_{k}$ of 
    $\vpol(E)$ is said to be \emph{irreducible} if the representation
    $(L_{-1},\ad_{\vert L_{0}})$ is irreducible.
\end{defi}
\begin{thm}\label{thm.main}
    Let $L=L_{-1}\oplus\cdots\oplus L_{k}$ be a graded subalgebra of 
    $\vpol(E)$. Then $L$ is maximal if and only if
    \begin{enumerate}
        \item  $L_{-1}=\vpol[-1](E)$,
    
        \item  $L$ is irreducible,
    
        \item  $L_{1}\neq 0$,
	\item  when $\K=\R$, the representation $(L_{-1},\ad_{\vert L_{0}})$ admits no complex 
	structure.
    \end{enumerate}
\end{thm}

The text below is organized as follows.
In Section~\ref{sec.cst}, we  prove the necessity of the first three  conditions above. 
Then, in Section~\ref{sec.mod}, we consider polynomial vector fields 
from a slightly modified point of view, in order to prove, in 
Section~\ref{sec.cplx}, the fourth condition given above.
We expose in Section~\ref{sec.IFFT} how the graded maximal subalgebras 
relate to the irreducible filtered Lie algebras of finite type, which 
were classified in~\cite{koba1}.
Using the classification of all irreducible infinite 
dimensional subalgebras of polynomial vector fields~(see for instance~\cite{guill, guill1, 
sing, koba3}, and references therein), we show in 
Section~\ref{sec.IFFTmax} that
all these algebras give rise to a canonical graded maximal subalgebra of 
polynomial vector fields.
\section{Constant vector fields and irreducibility}\label{sec.cst}
\begin{lem}
Let $L$ be a maximal subalgebra of $\vpol(E)$.  Then $\euler\in L$
if and only if $L$ is graded.
\end{lem}
\begin{proof}
    The sufficiency of the condition follows from the fact that 
    \[
    \K \euler + L
    \]
    is a Lie subalgebra when $L$ is graded.  In order to check the necessity 
of the condition, notice that 
    \[
    \ad(\euler)^k L\subset L, \quad \forall k\in\N,
    \] 
    gives a Vandermonde
    system allowing to compute the homogeneous components of a vector
    field $X\in L$.
\end{proof}
This proof is similar to the proof by Koecher (see~\cite[p.~354]{koe})
that any ideal of $\vpol(E)$ is graded.  We therefore state the
following remark.
\begin{rem}
    If $L$ is a subalgebra of $\vpol(E)$ that contains $\euler$, then 
any
    ideal of $L$ is graded.
\end{rem}

\begin{lem}\label{lem.prealg}
    Let $L_{-1}$, $L_{0}$ and $L_{+}$ be vector subspaces of
    $\vpol[-1](E)$, $\vpol[0](E)$ and $\bigoplus_{i\geq
    1}\vpol[i](E)$ respectively, such that
    \begin{enumerate}
        \item $L_{-1}\oplus L_{0}$ is a Lie subalgebra
    
	\item $[L_{-1},L_{+}]\subset L_{0}\oplus L_{+}$, and 
$[L_{0},L_{+}]\subset L_{+}$.
    \end{enumerate}
    Set $\descc[0](L_{+})=L_{+}$ and
    $\descc[k+1](L_{+})=[L_{+},\descc[k](L_{+})], (k\in \N)$.
    
    Then the smallest Lie subalgebra containing 
    \[
    L_{-1}\oplus L_0\oplus L_+
    \]
    is
    \[
    L_{-1}\oplus
    L_{0}\oplus\sum_{k\in\N}\descc[k](L_{+}).
    \]
    In particular, if $L_{+}\subset \vpol[1](E)$, the latter subalgebra is graded.
\end{lem}
\begin{proof}
    Using Jacobi identity, we check that 
$[L_{0},\descc[k](L_{+})]\subset
    \descc[k](L_{+})$ and consequently 
$[L_{-1},\descc[k](L_{+})]\subset
    \sum_{i=0}^k \descc[k](L_{+}) $ by induction on $k\geq 1$.  By 
definition,
    $[\descc[0](L_{+}),\descc[k](L_{+})]= \descc[k+1](L_{+})$ for all 
$k\in \N$.  Then, we check, by induction on $j\geq 0$, that 
$[\descc[j](L_{+}),
    \descc[k](L_{+})]\subset \descc[j+k+1](L_{+})$.

Therefore, $L_{-1}\oplus L_0\oplus\sum_{k\in\N} \descc[k](L_{+})$ is 
a Lie subalgebra. 
It is trivially the smallest one to contain the subspaces $L_{-1}$,
$L_{0}$ and $L_{+}$.
\end{proof}
\begin{defi}
    Let $F$ be a vector subspace of $\vpol[-1](E)$.  We set
    \[
    \n[i](F)=\{ X\in\vpol[i](E):\ad(F)^{i+1} X\subset F \}
    \]
    and
    \[
    \n(F) = \oplus_{i\geq -1} \n[i](F).
    \]
\end{defi}
Notice that $\n[-1](F)=F$ and that $\n[0](F)$
    is the intersection of the normalizer of $F$ and the subspace of
    linear vector fields.
\begin{prop}
    Let $L=\oplus_{i\geq -1} L_{i}$ be a graded subalgebra of
    $\vpol(E)$.  Then $\n(L_{-1})$ is an infinite dimensional graded
    subalgebra containing $L$.  Moreover, $\n(L_{-1})= \vpol(E)$ if and only if 
    $L_{-1}=\vpol[-1](E)$.
\end{prop}
\begin{proof}
    It is obvious that $[\n[i](L_{-1}),\n[j](L_{-1})]\subset
    \n[i+j](L_{-1})$.  Furthermore, if $L_{-1}=0$ or 
$L_{-1}=\vpol[-1](E)$,
    \[
    \n(L_{-1})=L_{-1}\oplus \bigoplus_{i\geq 0}\vpol[i](E).
    \]
    Now, if $h\in L_{-1}$, then, for every 
    polynomial function $p:E\rightarrow \K$, the field $x\mapsto 
    p(x)h$ belongs to $\n(L_{-1})$.
\end{proof}
\begin{cor}\label{cor.vect0}
    Let $L$ be a finite dimensional graded maximal subalgebra of
    $\vpol(E)$.  Then $\vpol[-1](E)\subset L$.
\end{cor}
\begin{cor}\label{cor.l1}
    Let $L$ be a finite dimensional graded maximal subalgebra of
    $\vpol(E)$.  Then $L_{1}\neq 0$.
\end{cor}
\begin{proof}
    Notice that $L$ cannot be made only of constant and linear vector 
    fields.  Indeed, it would then be included in the maximal 
    subalgebra~(\ref{sln}) presented in the introduction, for instance.
    Therefore, $L_{k}\neq 0$ for some $k>0$.  
    The conclusion follows from Corollary~\ref{cor.vect0}.
\end{proof}
\begin{prop}\label{prop.ideal0}
    Let $L$ be a finite dimensional graded maximal subalgebra of 
    $\vpol(E)$.  Then
    \[
    (L_{-1}=\vpol[-1](E), \ad_{\vert L_{0}})
    \]
    is an irreducible 
    representation of $L_0$.  
    It follows that any non trivial ideal of $L$ 
    contains every constant vector field.
\end{prop}

\begin{proof}[Proof of Proposition~\ref{prop.ideal0}]
  Let $F\neq\{0\}$ be a stable subspace of $L_{-1}$ under the action 
  of $L_0$.  
  
  The space
  \[
  L_{-1} \oplus \n[0](F) \oplus \bigoplus_{i\geq 1} 
  \{ X \in \vpol[i](E) : \ad(L_{-1})^i X\subset \n[0](F) \}
  \]
  satisfies the hypotheses of Lemma~\ref{lem.prealg}.  Its algebraic 
  closure is  an infinite dimensional proper subalgebra 
  containg $L$ properly, hence a contradiction.
  
  Let now $I$ be a non trivial ideal of $L$.  
  It contains at least one constant vector field since 
$[\vpol[-1](E),I]\subset I$.  
  It contains all of them since $I\cap L_{-1}$ is a stable subspace
  of $L_{-1}$.
\end{proof}

\section{A convenient model for polynomial vector 
fields}\label{sec.mod}
It will be useful to consider the spaces of multilinear symmetric 
mappings from $E\times \cdots \times E$ to $E$ instead of those of 
homogeneous polynomial vector fields.  We shall write
\[
\T[i](E)= S^{i+1}E^*\otimes E,\quad \text{ and }
\T(E)=\bigoplus_{i\geq -1} \T[i](E).
\]
To turn $\T(E)$ into a Lie algebra, we define as in~\cite{koba3} the 
following bracket operation.  If $t\in\T[p](E)$ and $t'\in\T[q](E)$ 
then $[t,t']\in\T[p+q](E)$ and
\begin{multline*}
[t,t'](x_{0},x_{1},\ldots,x_{p+q})= \\
\frac{1}{p!(q+1)!}\sum_{j} 
t(t'(x_{j_{0}},x_{j_{1}},\ldots,x_{j_{q}}),x_{j_{q+1}},\ldots,x_{j_{p+q}})\\
\mbox{-}\frac{1}{(p+1)!q!}\sum_{k} 
t'(t(x_{k_{0}},x_{k_{1}},\ldots,x_{k_{p}}),x_{k_{p+1}},\ldots,x_{k_{p+q}}) 
\end{multline*}
where both $j$ and $k$ run over all possible permutations of the 
$p+q+1$ first natural numbers.
\begin{prop}
    The map $T : \T(E)\rightarrow \vpol(E)$ defined 
 by
  \[
  T(M) : x\in E\mapsto -\frac{1}{(p+1)!} M(x,\ldots, x),\quad\forall 
  M\in \T[p](E)
  \]
  is an 
isomorphism of Lie algebras.
\end{prop}
\section{Absence of complex structure}\label{sec.cplx}
We now assume $\K=\R$ and prove, in Lemma~\ref{lem.vectCn}, the 
fourth condition of 
maximality of our main result.

 Let $E$ be a real vector space of even dimension and $J$ a complex 
 structure of $E$, i.e. an 
 endomorphism of $E$ such that $J^2=-\id$.  We denote by $E_{J}$ 
 the complex vector space defined by $E$ with the structure of 
 $\C$-module defined by
 \[
 (a+ i b) e := a e+b Je,\quad \forall a,b\in\R,\forall e\in E. 
 \]
 
 Define
 \begin{multline*}
 \T[p]^J(E)= \\ \{M\in\T[p](E)\vert
 J(M(x_{0},\ldots,x_{p}))=
 M(Jx_{0},x_{1}\ldots,x_{p}),
 \forall x_{0},\ldots,x_{p}\in E\}
 \end{multline*}
for all $p\geq -1$.  Then the subalgebra $\T^J(E)=\bigoplus_{i\geq 
-1}\T[p]^J(E)$ of $\T(E)$ is isomorphic to $\T(E_{J})$ as a real Lie 
algebra.  Indeed, the condition defining $\T^J(E)$ means that an 
application $M\in\T(E)$ is $\C$-multilinear on $E_{J}$.

\begin{lem}\label{lem.vectCn}
    Let $E$ be a real vector space and 
$L=\T[-1](E)\oplus\bigoplus_{j= 
    0}^k L_{j}$ a 
   graded subalgebra 
    of $\T(E)$.  Assume that $J$ is a complex structure of 
    $(L_{-1},\ad_{\vert L_{0}})$, i.e. 
    \begin{eqnarray*}
        &&[x_{0}, Jx_{-1}]=J[x_{0},x_{-1}],\quad \forall x_{0}\in 
        L_{0}, \forall x_{-1}\in L_{-1},\\
	&\text{and}& J^2 = -\id.
    \end{eqnarray*}
    Then
    \[
    L\subset \T^J(L_{-1})\subset \T(L_{-1})
    \]
    where both inclusions are strict.
\end{lem}
\begin{proof}
    Indeed, $L_{-1}=\T[-1]^J(E)=\T[-1](E)$.  The requirement for $J$ 
    to intertwine the action of $L_{0}$ on $L_{-1}$ precisely means 
    that $L_{0}\subset\T[0]^J(E)$.  If $L_{k-1}\subset \T[k-1]^J(E)$ 
    and $M\in L_{k}$, the equalities
    \begin{eqnarray*}
    J\circ M(x_{0},\ldots, x_{k}) &=& J( [M, x_{1}](x_{0},x_{2},\ldots, 
    x_{k})) \\
    &=& M(x_{1}, Jx_{0},x_{2},\ldots, x_{k})\\
    &=& M(Jx_{0},x_{1},\ldots, x_{k})
    \end{eqnarray*}
    show that $M\in\T[k]^J(E)$.
    
    The inclusions are strict because the dimension of $\T^J(L_{-1})$ 
    is infinite and because the dimension of $\T[p]^J(L_{-1})$, for 
    all $p\geq 0$, is  
    strictly less than that of $\T[p](L_{-1})$.
\end{proof}

This lemma generalizes the construction used in~\cite{bolec} to show 
that a subalgebra of infinitesimal conformal transformations, 
isomorphic to $\alg{so}{3,1,\R}$,  
is not maximal in $\vpol(\R^2)$.
\section{Irreducible filtered algebras of finite type}\label{sec.IFFT}
Let $L=\oplus_{j=-1}^k L_{j}$ be a graded maximal subalgebra of 
polynomial vector fields.
In the last section, we have shown that $L$ possesses interesting 
properties.  It actually belongs to a broader class of Lie algebras 
studied in~(\cite[Theorem~1, p.~875]{koba1}).  

This theorem describes the structure of some filtered finite 
dimensional Lie algebras together with a group of automorphisms.

We shall only associate the trivial group $\{\id\}$ to such an algebra. 
Furthermore, the reader may find worth noticing that the algebras we 
consider carry the filtration which is naturally associated to their grading 
and that the other hypotheses of the theorem are satisfied in view of 
the first three conditions required in our main result for a 
subalgebra to be maximal.

For the sake of simplicity, we shall name algebras described by this theorem 
\emph{Irreducible filtered algebras of finite type}, as it was done 
in~\cite{koba3}, or simply write \emph{IFFT-algebras}.

As a consequence of the mentioned result, we know that $L$ is simple 
and is \emph{of order two}, i.e. $L=L_{-1}\oplus L_{0}\oplus L_{1}$.
Moreover, there exists a unique element $e\in L$ such that $L_{p}$ is 
the eigenspace of $\ad(e)$ associated to the eigenvalue $p$.  This 
element is in the center of $L_{0}$.  We shall name it the \emph{Euler} 
element of $L$.
Finally, $L_{-1}$ and $L_{1}$ are dual to each other as modules of 
$L_{0}$ with respect to the Killing form of $L$.

On the one hand, Kobayashi and Nagano gave a list of the admissible algebras 
and detailed in each case the associated grading.  The pairs 
 $(L,e)$ where 
$L$ is a real IFFT-algebra and $e$ its Euler element 
are classified in~\cite[pp.~892--895]{koba1}.
On the other hand, to any graded algebra $L=\bigoplus_{k\geq -1} L_{k}$,
they associated in a natural way a graded subalgebra of $\T(L_{-1})$~%
(see~\cite[p.~683]{koba3}). The reader may compare this construction 
with that of Gradl~(\cite{gradl}).
In the case of $L=L_{-1}\oplus L_{0}\oplus L_{1}$, this is done by
the following monomorphism 
$\phi:L\rightarrow\T(L_{-1})$ :
\[
\left\{
\begin{array}{lcl}
    \phi_{\vert L_{-1}}&=&\id\\
    \phi_{\vert L_{0}} &=&\ad_{\vert L_{0}}\\
    \phi(M)&=&(x,y)\mapsto [[M,x],y],\quad\forall M\in L_{1},\forall 
    x,y\in L_{-1}.
\end{array}
\right.
\]
Notice that this is the only way to proceed provided 
the value of $\phi$ on $L_{-1}$ is set to $\id$.
\section{IFFT-algebras are maximal}\label{sec.IFFTmax}
In this section, we prove the sufficiency of the conditions 
given in our main result for a subalgebra of polynomial vector fields 
to be maximal.  

We first assume that $E$ is a complex vector space and 
$L=L_{-1}\oplus L_{0}\oplus L_{1}$ an irreducible graded subalgebra 
of $\T(E)$ such that $L_{-1}=\T[-1](E)$ and $L_{1}\neq 0$.
Then we shall show how the proof adapts to the real case.

Let $L'$ be a subalgebra of $\T(E)$ such that $L'\supset L$. Then 
$L'$ 
is graded and irreducible, since $L$ is.

If $L'$ is finite dimensional, one sees,  by 
using the description of the IFFT-algebras~(see 
Section~\ref{sec.IFFT}), that $L'_{1}=L_{1}$, that $L'$ is 
simple and eventually that $L=L'$, since $[L'_{-1},L'_{1}]=L_{0}$.

Therefore, if $L'$ contains properly $L$, then it must be infinite 
dimensional.  It possesses two additional properties, consequences of 
the following result.
\begin{prop}[\protect{\cite[p.~688]{koba3}}]
    Let $\bigoplus_{p\geq -1} G_{p}$ be an irreducible graded Lie 
    algebra of infinite type or finite type of order $\geq 2$ over a 
    field of characteristic 0.  Then $G_{0}$ is reductive and 
    $[G_{-1},G_{1}]$ contains the semi-simple part of $G_{0}$.
\end{prop}
\begin{itemize}
    \item  $L'_{0}$ is reductive and has a non trivial center (the 
multiples of the 
    identity transformation of $L_{-1}$).  

    \item $L'$ is still simple.  Indeed, if $I$ is an ideal of $L'$ 
    then $I\supset L$ (see Proposition~\ref{prop.ideal0}), which 
    implies that $I$ contains the multiples of the identity 
    transformation of $L_{-1}$ and in turn that $I\supset L'_{j}$ for 
    all $j\neq 0$.  Since $[L'_{-1},L'_{1}]$ contains the semi-simple 
    part of $L'_{0}$, the conclusion follows. 
\end{itemize}

In order to prove the maximality of $L$, the remaining point is to ensure that 
$L'=\T(E)$.  The key result is due to Cartan.  We refer the reader to 
the works of Guillemin, Quillen, Singer, Sternberg, Kobayashi and 
Nagano~(\cite{guill, guill1, sing, koba3}).

This result states that the only irreducible infinite dimensional 
graded 
subalgebras of $\vpol(E)$ are
\begin{enumerate}
    \item  $\vpol(E)$ itself,

    \item  the divergence-free vector fields,

    \item  the Hamiltonian vector fields with respect to a symplectic 
    form given on $E$, provided $E$ is even 
    dimensional,

    \item  the last two subalgebras supplemented with the multiples 
of 
    the Euler vector field.
\end{enumerate}

But the subalgebras described in~(4) are not simple, and those in~(2) 
and~(3) have a simple linear part.

Hence the proof.

Now, when $E$ is a real vector space and $L$ and $L'$ as above, one 
proceeds in the same way to prove the simplicity of $L'$, noticing 
that both $L_{0}$ and $L'_{0}$ still have a one dimensional center.  

Indeed, 
if $x_{0}$ is central in one of these two subalgebras, then 
$\ad(x_{0})$ intertwines the action of $L_{0}$ on $L_{-1}$.  Since 
the representation $(L_{-1},\ad_{\vert L_{0}})$ admits no complex 
structure, Schur's lemma ensures that $ad(x_{0})_{\vert L_{-1}}$ is a multiple of the 
identity transformation of $L_{-1}$.  Therefore, $\dim Z(L_{0})=1$.

The description of irreducible infinite dimensional graded 
subalgebras 
of $\vpol(E)$ is essentially due to Matsushima.  It can be found in~\cite{koba3, matsu}.

Two cases arise whether $L'_{0}\otimes \C$ acts 
irreducibly on $E\otimes\C$ or not.  In the first case, 
$L'$ should be one of the real analogues of the Cartan algebras 
listed 
above.  But in the second, $E$ admits a complex 
structure as a $L_{0}$ module, which contradicts the hypotheses.

Theorem~\ref{thm.main} is proved.

In order to complete our search for maximal subalgebras of polynomial vector 
fields over a given \emph{real} vector space, we need to be able to 
identify in the tables given in~\cite[pp.~892--895]{koba1},
the algebras such that 
the representation $(L_{-1},\ad_{\vert L_{0}})$ admits a 
complex structure.

\begin{prop}
    Let $L_{-1}$ be a real vector space and $L=L_{-1}\oplus L_{0}\oplus 
L_{1}$  
    an IFFT-algebra.   
    Then $(L_{-1},\ad_{\vert L_{0}})$ admits a complex structure
    if and only if the algebra $L$ admits a complex structure.
\end{prop}
\begin{proof}
    The sufficiency of the condition is obvious.  Notice that a 
    complex structure on $L$ stabilizes the eigenspaces of $e$.
    
    Let $J_{-1}$ be a complex structure on $(L_{-1},\ad_{\vert 
    L_{0}})$.  Let $J_{1}:L_{1}\rightarrow L_{1}$ be the adjoint of 
    $J_{-1}$ with respect to the $Killing$ form $\beta$ of $L$, i.e.
    \[
    \beta(J_{1} x_{1}, x_{-1})=\beta(x_{1}, J_{-1} x_{-1}), \quad 
    \forall x_{-1}\in L_{-1}, \forall x_{1}\in L_{1}.
    \]
    The so defined $J_{1}$ intertwines the action of $L_{0}$ on 
$L_{1}$.
    Moreover,
    \[
    [J_{1} x_{1}, x_{-1}]=[ x_{1}, J_{-1} x_{-1}], \quad \forall 
    x_{-1}\in L_{-1}, \forall x_{1}\in L_{1}.
    \]
    Indeed, for all $x_{-1}, y_{-1}\in L_{-1}$ and $x_{1}, y_{1}\in 
    L_{1}$,
    \begin{eqnarray*}
        \beta([[x_{1},J_{-1} x_{-1}], y_{-1}], y_{1})&=& 
        \beta(J_{-1} x_{-1},[x_{1},[y_{1},y_{-1}])\\
	&=&\beta (x_{-1}, [J_{1} x_{1}, [y_{1},y_{-1}]])\\
	&=&\beta([[J_{1} x_{1}, x_{-1}],y_{-1}],y_{1}).
    \end{eqnarray*}
    Define 
    \[
    J_{0}:L_{0}\rightarrow \T[0](L_{-1}):A\mapsto  A\circ J_{-1}.
    \]
    This map is actually valued in $L_{0}$ since
    \begin{eqnarray*}
        (J_{0} [x_{1},x_{-1}])y_{-1}&=& [[x_{1},x_{-1}], 
        J_{-1}y_{-1}] \\
	&=& [[x_{1}, J_{-1}y_{-1}],x_{-1}]\\
	&=& [[J_{1}x_{1},y_{-1}],x_{-1}]\\
	&=& [J_{1}x_{1},x_{-1}]y_{-1}
    \end{eqnarray*}
    for all $y_{-1}\in L_{-1}$.
    
    The map $J:L\rightarrow L$ defined by its restrictions $J_{i}$ to 
    $L_{i}$ $(i=-1,\ldots,1)$ is then a complex structure of 
    $L$ as a Lie algebra.
\end{proof}

The statement ``IFFT-algebras are maximal'' should be 
taken in the following sense. In the tables given in~\cite{koba1}, 
 one can distinguish complex algebras  from 
real ones  admitting no complex structure.
The latter give rise to maximal subalgebras of the real algebra 
$\T(L_{-1})$.
One may consider the former as Lie subalgebras of the real Lie 
algebra $\T(L_{-1})$, in which case Lemma~\ref{lem.vectCn} shows that 
they are not maximal.
They are maximal when regarded in their natural position of complex 
subalgebras of the 
complex Lie algebra $\T(L_{-1})$.

\section*{Acknowledgements}
We are very grateful to C. Duval, M. De Wilde, P. Lecomte  and V. 
Ovsienko for their 
interest and suggestions.

The first author thanks the  Belgian 
National Fund for Scientific Research (FNRS) for his Research 
Fellowship.

Institute of Mathematics, B37
\\University of Li\`ege,
\\B-4000 Sart Tilman
\\Belgium
\\
\\mailto:f.boniver@ulg.ac.be
\\mailto:p.mathonet@ulg.ac.be

\begin{thebibliography}{100}

\bibitem{dlo} C. Duval, P. Lecomte and V. Ovsienko, 
``Conformally equivariant quantization: existence and uniqueness,'' 
Ann. Inst. Fourier (Grenoble)  49  no. 6, 1999--2029 (1999)
\bibitem{lo} P. B. A. Lecomte and V. Ovsienko, 
``Projectively equivariant symbol calculus,'' 
Lett. Math. Phys. 49 no. 3, 173--196  (1999) 
\bibitem{bolec} F.  Boniver and P.  B.  A.  Lecomte, 
``A remark about the Lie algebra of infinitesimal conformal transformations of the Euclidian space,'' 
Bull.  London Math.  Soc. 32, 263--266 (2000)
\bibitem{kan} I. L. Kantor,
``Classification of irreducible transitively differential groups,''
Sov. Math. Dokl. 5, 1404--1407 (1964)
\bibitem{post} G.  Post, 
``A class of graded Lie algebras of vector fields and first order differential 
operators,''
J. Math. Phys. 35, no. 12, 6838--6856 (1994)
\bibitem{koba3} S. Kobayashi and T. Nagano, ``On Filtered Lie Algebras 
and Geometric Structures III,'' J. Math. Mech. 14 no. 4, 
679--706 (1965)
\bibitem{koba1} S. Kobayashi and T. Nagano, ``On Filtered Lie Algebras 
and Geometric Structures I,'' J. Math. Mech. 13 no. 5, 
875--907 (1964)
\bibitem{guill} V. W. Guillemin, D. Quillen, S. Sternberg, ``The 
classification of the irreducible complex algebras of infinite type,'' J. 
Analyse Math. 18, 107--112 (1967)
\bibitem{guill1} V. W. Guillemin and S. Sternberg, ``An algebraic model 
of transitive differential geometry,'' Bull. Amer. Math. Soc., vol. 
70, 16--47 (1970)
\bibitem{sing} I. M. Singer, S. Sternberg, ``The infinite groups of 
Lie and Cartan.  I. The transitive groups,'' J. Analyse Math. 15, 
1--114 (1965)
\bibitem{koe} M. Koecher, 
``Gruppen und Lie-Algebren von rationalen Funktionen,''
Math.  Z.  109, 349--392 (1969)
\bibitem{gradl} H. Gradl, ``Realization of semisimple Lie algebras with 
polynomial and rational vector fields,'' Comm. Algebra 21 no. 11, 
4065--4081 (1993)
\bibitem{matsu} Y. Matsushima, ``Sur les alg\`ebres de Lie 
semi-involutives,'' Colloque de topologie de Strasbourg, 
1--17 (1954)


\end{thebibliography}
\end{document}